\documentclass[1p]{elsarticle}

\usepackage{lineno,hyperref}
\modulolinenumbers[5]
\usepackage{amsmath,amssymb}
\usepackage{algorithm}
\usepackage[noend]{algorithmic}
\usepackage{float}
\usepackage {tikz}

\newtheorem{theorem}{Theorem}
\newtheorem{lemma}{Lemma}
\newtheorem{corollary}{Corollary}

\newproof{pf}{Proof}
\newcommand{\uxi}{\ensuremath{\underline{\xi}}}

\begin{document}

\begin{frontmatter}

  \title{Robust Single Machine Makespan Scheduling with Release Date Uncertainty}

  \author[mymainaddress]{Oliver Bachtler}
  \ead[url]{https://www.mathematik.uni-kl.de/opt/personen/mitglieder/}
  \ead{bachtler@mathematik.uni-kl.de}

  \author[mymainaddress]{Sven O. Krumke}
  \ead[url]{https://www.mathematik.uni-kl.de/opt/personen/leitung/krumke/}
  \ead{krumke@mathematik.uni-kl.de}
  
  \author[mymainaddress]{Huy Minh Le\corref{mycorrespondingauthor}}
  \cortext[mycorrespondingauthor] {Corresponding author}
  \ead[url]{https://www.mathematik.uni-kl.de/opt/personen/mitglieder/}
  \ead{leminh@mathematik.uni-kl.de}

\address[mymainaddress]{Department of Mathematics, Technische Universit\"at Kaiserslautern, Paul-Ehrlich-Str. 14, 67663~Kaiserslautern, Germany}

\begin{abstract}
  This paper addresses the robust single machine makespan scheduling with uncertain release dates of the jobs.  The release dates take values within know intervals. We use the concept of Gamma-robustness in two different settings and address both 
 the robust absolute and robust regret criteria. Our main results are polynomial time algorithms which have the same running time ($O(n \log n)$) as the best algorithms for the non-robust case.
\end{abstract}


\end{frontmatter}

\section{Introduction}
Scheduling theory is an important area of operations research with a wealth of applications, e.g., in management, production, computer systems, construction, etc. Typically scheduling theory deals with the task that a finite set of jobs needs to be processed by a system with limited resources. One of  the most popular objectives that has been used in scheduling theory is the makespan objective. We refer to the books by Brucker~\cite{Brucker07}, Pinedo~\cite{Pinedo08}, and the survey papers of Lenstra et.\ al.~\cite{Lenstra77} and~\cite{Lenstra93} for references.

In practice decision makers are interested in hedging against the worst possible scenarios.  A solution is often compared with an optimal solution that could have been obtained if the actual realization of the uncertain parameters had been available. In order to meet these requirements, the robust optimization framework has been proposed in~\cite{Bental09,Kouvelis97,Soyster73}. Robust optimization looks for solutions in a context where one is faced with imprecise, uncertain and generally incompletely known parameters of a problem. In most robust optimization models, we compute a solution minimizing the largest cost (absolute criterion) or the largest deviation from the optimum (regret criterion) where the maximum is taken over all scenarios which describe the uncertainty set.  Both criteria are well-known in the theory of decision making under uncertainty, where no probability distribution is at hand.

In this paper, we consider scheduling jobs on a single machine. Processing times are precisely known but all release dates of the jobs are uncertain. A class of scheduling problems with uncertain parameters is discussed by Kasperski and Zielinski~\cite{Kasperski14}, where the complexity of various (regret) scheduling problems and some algorithms for solving them are described for both the absolute and regret criteria. The budgeted uncertainty for some scheduling problems is considered in~\cite{Bougeret19,Bruni18,DanielsKouvelis95,Luetal12}. The concept of Gamma-robustness for combined scheduling-location problems is studied in~\cite{Krumke20} in the sense that the total deviation of the uncertainty parameters cannot exceed some threshold. For other approaches to robust optimization, we refer to~\cite{Bental09,Soyster73}.

The paper is organized as follows. Section~\ref{section_notation_framework} introduces notation and the general framework. In Section~\ref{section_Robustabsolutecriterion}, we use the robust absolute approach to solve the scheduling problem with uncertain parameters in two different settings.  We provide $O(n\log n)$ time algorithms for both cases.  Section~\ref{section_Robustregretcriterion} is dedicated to the robust regret approach. We first give an $O(n^2)$ time algorithm to find an optimal schedule for this problem and then describe an implementation with an improved running time of~$O(n\log n)$. 

\section{Notation and the General Framework}
\label{section_notation_framework}

The single machine makespan scheduling problem ($1|r_j|C_{\max}$) is defined as follows: A set of jobs $\mathcal{J} = \{ 1, \dots, n \}$ needs to be scheduled on a single machine.  Each job $j\in \mathcal{J}$ has a processing time $p_j > 0$ and a release date $r_j >0$ before which it cannot be scheduled.  The jobs need to be processed without preemption and the machine can handle only one job at any point in time.  Let $\Pi$ denote the set of all permutations $\pi$ of the jobs where $\pi(i) = j$ states that job $j$ is placed at position $i$. The goal is to find an optimal schedule among all permutations of jobs such that the makespan value is minimized, i.e., a permutation~$\pi$ minimizing~$C_{\max} (\pi)$, where $C_{\max} (\pi) = C_{\pi(n)} (\pi)$ with $C_{\pi(i)} (\pi)$ is the completion time of job~$\pi(i)$. Completion times of the jobs are computed as follows: $C_{\pi(1)} (\pi) = r_{\pi(1)} + p_{\pi(1)}$ and for all $i = 2, \dots, n$,
\begin{align*}
  C_{\pi(i)} (\pi) = \max \{C_{\pi(i-1)} (\pi), r_{\pi(i)} \} + p_{\pi(i)}.
\end{align*}
Lenstra et al. \cite{Lenstra77} showed that an optimal schedule can be found by applying the Earliest Release Date Rule (ERD).  Thus, by sorting the jobs in non-decreasing order of their release dates we obtain an algorithm with running time~$O(n\log n)$.

The robust version in this paper is obtained by modifying certain parts of the input data to be uncertain: The processing times $p_j > 0$ are known, but the release dates $r_j$ are uncertain and may take any value within known intervals $[\underline{r}_j, \overline{r}_j]$. We call this problem the \emph{robust single machine makespan scheduling problem} (briefly, Robust $1|r_j|C_{\max}$) with uncertainty in the release dates.  We handle the robust version by using the concept of Gamma-robustness in two different settings of the uncertainty set. These are the continuous and discrete versions of budgeted uncertainty: In the first setting the total deviation of the uncertain parameters cannot exceed a given threshold~$\Gamma$, which is not necessarily integer, i.e.,
\begin{align*}
  U_1 = \bigl\{ \, \xi \mathrel{\vert} r_j^{\xi} \in [\underline{r}_j, \overline{r}_j]  \text{ for all } j \in \mathcal{J}  \text{ and } \sum_{j \in \mathcal{J}} (r_j^{\xi} - \underline{r}_j) \leq \Gamma \,\bigr\}.
\end{align*}
Alternatively, we regard the setting in which the number of jobs that can deviate from their corresponding lower bounds does not exceed some threshold $\Gamma$, in which $\Gamma$ is considered to be an integer and at least~$1$, i.e.,
\begin{align*}
  U_2 = \bigl\{ \, \xi \mathrel{\vert} r_j^{\xi} \in [\underline{r}_j, \overline{r}_j] \text{ for all }  j \in \mathcal{J} \text{ and } |\{ j : r_j^{\xi} \neq \underline{r}_j \} | \leq \Gamma \,\bigr\}.
\end{align*}
Each realization $\xi$ in the uncertainty sets is called a \emph{scenario}.

Let $\pi$ be a permutation (i.e.\ a schedule).  For each job~$j$, let $r_j^{\xi}$ and $C_j(\pi, \xi)$ be the release date and completion time of job~$j$ in~$\pi$ under scenario~$\xi$, respectively. Let $\overline{\xi}$ be the (hypothetical) scenario in which all release dates take their latest values, i.e., $r_j^{\overline{\xi}} = \overline{r}_j$. Note that it may be the case that $\bar{\xi}$ is not in $U_1$ or $U_2$, since $\Gamma$ may not be large enough. However, the (hypothetical) scenario $\bar{\xi}$ will be useful in the later analysis. A job $j \in \mathcal{J}$ is called a \textit{critical job} for a permutation~$\pi$ under scenario~$\xi$ if $r_j^{\xi}+ \sum_{s=i}^{n} p_{\pi(s)} = C_{\max} (\pi, \xi)$, where $\pi(i) = j$.  
Note that such a critical job always exists: Let $j$ be the last job in~$\pi$ satisfying $C_j(\pi,\xi)=r_j^\xi+p_j$ (which exists since the first job always satisfies this condition).  The jobs after~$j$ are then processed without idle time and, thus, $j$ is critical.

In the next two sections we consider the robust absolute and robust regret versions of our scheduling problem.
All results obtained in these sections apply to both uncertainty sets $U_1$ and $U_2$, and the corresponding proofs are similar.
To illustrate both cases, we deal with $U_1$ in Section~\ref{section_Robustabsolutecriterion} and derive our results in Section~\ref{section_Robustregretcriterion} for $U_2$.

\section{Robust Absolute Criterion}
\label{section_Robustabsolutecriterion}

Consider the uncertainty set~$U_1$. The goal of the robust absolute single machine makespan scheduling problem (Robust Absolute $1|r_j|C_{\max}$) is to find a schedule such that its cost in the worst-case scenario is minimized, i.e., we wish to find a permutation~$\pi$ minimizing $\max_{\xi \in U_1} C_{\max} (\pi, \xi)$.  Without loss of generality, we may assume that $\Gamma \geq \max_{j \in \mathcal{J}} \{ \overline{r}_j - \underline{r}_j \}$ since otherwise one can reduce the uncertain intervals of release date to $[\underline{r}_j, \underline{r}_j + \Gamma]$ for any job~$j$ that satisfies $\overline{r}_j - \underline{r}_j > \Gamma$.

\begin{lemma}
	\label{worst_case_scenario_Gamma}
	For a fixed permutation $\pi$, let $j^{*}$ be a critical job for $\pi$ under the hypothetical scenario $\overline{\xi}$. There exists a worst-case scenario $\xi^{*}\in U_1$ for $\pi$ such that under this scenario the following holds: 
	\begin{enumerate}[(a)]
        \item The job~$j^{*}$ is released as late as possible and the release dates of all other jobs take the corresponding lower bounds, i.e., $r_{j^{*}}^{\xi^{*}} = \bar{r}_{j^{*}}$ and $r_{j}^{\xi^{*}} = \underline{r}_j$, for $j \neq j^{*}$;\label{upper_bound_worst_case_scenario_Gamma}
		\item The job~$j^{*}$ is a critical job for~$\pi$ under scenario $\xi^{*}$. \label{critical_job_worst_case_scenario_Gamma}
	\end{enumerate}
      \end{lemma}

      \begin{pf}
        Let $\xi$ be a worst-case scenario for~$\pi$.  Since increasing the release dates can never decrease the makespan, it follows that $C_{\max}(\pi,\xi) \leq C_{\max}(\pi,\bar{\xi})$. Since $j^*$ is critical for~$\pi$ under scenario~$\bar{\xi}$, it follows that there is no idle time in $\pi$ once $j^*$ is started, and the makespan $C_{\max}(\pi,\bar{\xi})$ depends solely on the release date $\bar{r}_{j^*}$ and the processing times of the jobs processed in$\pi$ after and including job~$j^*$. Consequently, the makespan does not decrease if we reduce the release dates of all jobs $j\neq j^*$ to the lowest possible value~$\underline{r}_j$. This yields a valid scenario $\xi^*$ which has the same makespan as the hypothetical scenario $\bar{\xi}$, and which thus must be a worst-case scenario for $\pi$. It remains to show that $j^*$ is actually a critical job for~$\pi$ under scenario~$\xi^*$. According to the construction of scenario~$\xi^{*}$, we can see that $r_{j^{*}}^{\xi^{*}} = \overline{r}_{j^{*}}$ and $r_j^{\xi^{*}} \leq \overline{r}_j$, for all $j \neq j^{*}$. Let $\pi(i^{*}) = j^{*}$ and $\pi(i) = j$ where $j\neq j^{*}$. Since $j^{*}$ is a critical job for~$\pi$ under scenario~$\bar{\xi}$, we obtain $r_{\pi(i^{*})}^{\xi^{*}} + \sum_{s = i^{*}}^{n} p_{\pi(s)} = \bar{r}_{\pi(i^{*})} + \sum_{s = i^{*}}^{n} p_{\pi(s)} \geq \bar{r}_{\pi(i)} + \sum_{s = i}^{n} p_{\pi(s)} \geq r_{\pi(i)}^{\xi^{*}} + \sum_{s = i}^{n} p_{\pi(s)} $ for all $i \neq i^*$. This implies that $j^{*}$ is a critical job for $\pi$ under scenario $\xi^{*}$. \qed 
      \end{pf}

      \begin{theorem}\label{thm:1}
	An optimal schedule~$\pi$ for Robust Absolute $1|r_j|C_{\max}$ is obtained by sorting all jobs in non-decreasing order of latest  possible release dates, i.e., 
        \begin{align*}
          \bar{r}_{\pi (1)} \leq \bar{r}_{\pi (2)} \leq \cdots \leq \bar{r}_{\pi (n)} .
        \end{align*}
      \end{theorem}

\begin{pf}
  Let $\pi$ be a schedule as in the claim and $j^{*}$ with $j^{*} = \pi(i^{*})$ be a critical job for~$\pi$ under scenario~$\bar{\xi}$. If $\xi^*$ is the scenario described in Lemma~\ref{worst_case_scenario_Gamma}, then $\xi^*$ is a worst-case scenario for $\pi$ and $j^{*}$ is a critical job for~$\pi$ under $\xi^*$. If $\pi' \neq \pi$ is another schedule, then by the Pigeonhole Principle there exists a job~$j$ whose position in~$\pi$ is in~$[i^{*}, n]$, and its position in~$\pi'$ is in~$[1, i^{*}]$. We choose the job~$j$ with smallest index in~$\pi'$ satisfying this property. Let $\pi(i) = j$ and $\pi'(i') = j$.  Then $i' \leq i^{*}\leq i$. Thus,
  \begin{align*}
    \max_{\xi \in U_1} C_{\max} (\pi', \xi) \geq \bar{r}_j + \sum_{s=i'}^{n} p_{\pi'(s)} \geq \bar{r}_{j^{*}} + \sum_{s=i^{*}}^{n} p_{\pi(s)} = \max_{\xi \in U_1} C_{\max} (\pi, \xi).    
  \end{align*}
  The second inequality holds since the permutation~$\pi$ is sorted in non-decreasing order of latest release dates and $j$ is chosen as the smallest index in~$\pi'$. \qed
\end{pf}

Theorem~\ref{thm:1} implies the polynomial time solvability of Robust Absolute $1|r_j|C_{\max}$ with uncertainty set~$U_1$.  Moreover, it is straight-forward to see that statement of Lemma~\ref{worst_case_scenario_Gamma} still holds for the uncertainty set~$U_2$.  Thus, Theorem~\ref{thm:1}  carries over to~$U_2$ as well. We obtain the following result:
\begin{corollary}
  The problem Robust Absolute $1|r_j|C_{\max}$ with uncertainty set~$U_1$ or~$U_2$ can be solved optimally in time~$O(n\log n)$.
\end{corollary}

\section{Robust Regret Criterion}
\label{section_Robustregretcriterion}

We consider the minimum regret version of~$1|r_j|C_{\max}$ and we start with uncertainty set~$U_2$.  The goal of Robust Regret $1|r_j|C_{\max}$ is to find an optimal regret schedule such that the makespan regret (opportunity loss) overall scenarios is minimal, i.e., to find a permutation~$\pi$ minimizing
\begin{align*}
R(\pi) = \max_{\xi \in U_2} R(\pi, \xi)= \max_{\xi \in U_2} (C_{\max} (\pi, \xi) - \text{M}(\xi)),
\end{align*}
where $\text{M}(\xi)=\min\{\, C_{\max}(\pi, \xi) : \pi\in\Pi\,\}$ is the optimal makespan value for the problem in a known scenario $\xi \in U_2$.  Recall that for a fixed scenario $\xi$ the problem can be solved in time $O(n \log n)$ by applying the ERD rule.  A job $j$, with $\pi(i) = j$, is called a \textit{critical-regret} job for $\pi$ if it is a critical job in a worst-case regret scenario $\xi$, i.e., $R(\pi) = r_j^{\xi} + \sum_{s = i}^{n} p_{\pi(s)} - \text{M}(\xi)$.

\begin{lemma}
  \label{increaseasufficientlysmallincrease}
  Fix a permutation~$\pi$ and a scenario~$\xi\in U_2$.  
  If $\xi'$ is the (hypothetical) scenario obtained from~$\xi$ by increasing the release date of a job~$j$ by some value $\varepsilon\geq 0$, then $C_{\max}(\pi,\xi')\leq C_{\max}(\pi,\xi)+\varepsilon$.
  In particular, $C_{\max}(\pi,\xi')= C_{\max}(\pi,\xi)+\varepsilon$ if $j$ is critical for $\pi$ under $\xi$ and $M(\xi')\leq M(\xi)+\varepsilon$.
\end{lemma}

\begin{pf}
	The inequality in the claim follows directly by an induction on the recursive definition of the makespan, since the completion time of any job is delayed by at most~$\varepsilon$.
	The claimed equality, if job~$j$ is critical for $\pi$ under~$\xi$, holds as it remains critical under scenario~$\xi'$ and, hence, the completion time of any job from~$j$ onward in~$\pi$ is delayed by exactly~$\varepsilon$.
	For the statement regarding the makespan, let $\pi^*$ be an optimal permutation for~$\xi$.
	The makespan of~$\pi^*$ increases by at most~$\varepsilon$ in the transition to~$\xi'$ and the optimum permutation for~$\xi'$ cannot have larger makespan.
	\qed
\end{pf}

Let $\xi_{j}$ be the scenario associated with job~$j$ such that:
\begin{enumerate}[(i)]
\item The release date of job $j$ takes the corresponding upper bound. \label{condition_upperbound}
\item The other release dates of jobs take the corresponding lower bounds. \label{condition_lowerbound}
\end{enumerate}
Such a scenario was already used in Lemma~\ref{worst_case_scenario_Gamma} in the special case that $j$ is a critical job.  Clearly, the above constructed scenario satisfies $\xi_j \in U_2$ for all $j \in \mathcal{J}$.

\begin{lemma}
  \label{worstcaseregretscenario}
  Let $\pi$ be a fixed permutation of the jobs.  Let $j\in \mathcal{J}$ be the job maximizing~$R(\pi, \xi_j)$.  Then $\xi_{j}$ is a worst-case regret scenario for $\pi$ in $U_2$.
\end{lemma}

\begin{pf}
  Let $\xi\in U_2$ be a worst-case regret scenario for $\pi$ and $j^{*}$ be a critical-regret job for $\pi$ under scenario~$\xi$.  By increasing the release date of job~$j^*$ by  $\varepsilon:=\bar{r}_{j^{*}} -r_{j^{*}}^{\xi}\geq 0$, we obtain a new scenario~$\xi'$ whose regret value does not decrease, i.e., $R(\pi,\xi)\leq R(\pi,\xi')$.  To see this, apply Lemma~\ref{increaseasufficientlysmallincrease} which yields that the makespan of~$\pi$ increases by exactly~$\varepsilon$ (since $j^*$ is a critical-regret job) and the optimum makespan increases by at most~$\varepsilon$.

  Scenario $\xi'$ may not be contained in~$U_2$.  Now let us show that $\xi'$ can be transformed into the valid scenario $\xi_{j^{*}}\in U_2$ without decreasing the regret. Consider any job $j' \in \mathcal{J}$ that differs from job~$j^{*}$. If we decrease the release date of this job, then the value $C_{\max}(\pi, \xi')$ does not change.  Moreover, the value $M(\xi')$ does not increase, hence the regret value $R(\pi,\xi')$ does not decrease. Thus, if we replace the release dates of all jobs~$j'\neq j^*$ by the corresponding lower bounds~$\underline{r}_{j'}$, the scenario $\xi' = \xi_{j^{*}}$ still remains a worst-case regret scenario for~$\pi$. By the assumption that~$j$ maximizes, $R(\pi,\xi_j)$, i.e., $R(\pi, \xi_j) \geq R(\pi, \xi_{j'})$, for all $j' \in \mathcal{J}$, we can conclude that $\xi_{j}$ is a worst-case regret scenario for~$\pi$.
  \qed
 \end{pf}

\begin{corollary}
	\label{finiteworstcaseregretscenario}
	The set $\{ \, \xi_{j} \; | \; j \in \mathcal{J} \,\} \subseteq U_2$ contains a worst-case regret scenario for every permutation $\pi$.
      \end{corollary}
      \begin{pf}
        For any permutation $\pi$, Lemma \ref{worstcaseregretscenario} shows that the scenario $\xi_j$ with $R(\pi, \xi_j)$ maximum is a worst-case regret scenario.\qed
      \end{pf}

\begin{theorem}
  \label{alg_regret}
  An optimal schedule~$\pi$ for Robust Regret $1|r_j|C_{\max}$ is obtained by sorting all jobs in non-decreasing order of $\bar{r}_j - \text{M}(\xi_j)$, i.e.,
  \begin{align*}
    \bar{r}_{\pi (1)} - \text{M}(\xi_{\pi(1)}) \leq \bar{r}_{\pi (2)} - \text{M}(\xi_{\pi(2)}) \leq \cdots \leq \bar{r}_{\pi (n)} - \text{M}(\xi_{\pi(n)}) .
  \end{align*}
\end{theorem}

\begin{pf}
  Let $\pi$ be a permutation as in the claim and let $j^{*}$ with $j^{*} = \pi(i^{*})$ be a critical-regret job for~$\pi$ under the worst-case regret scenario~$\xi_{j^{*}}$. If $\pi' \neq \pi$ is another permutation, then there exists a job~$j$ whose position in~$\pi$ is in $[i^{*}, n]$, and its position in~$\pi'$ is in $[1, i^{*}]$. We choose a job~$j$ with smallest index in~$\pi'$ satisfying this property. Let $\pi(i) = j$ and $\pi'(i') = j$.  Thus, $i' \leq i^{*}\leq i$ and
  \begin{align*}
    R(\pi') \geq R(\pi', \xi_{j}) &= C_{\max} (\pi', \xi_j) - \text{M}(\xi_{j})\\
            &\geq \bar{r}_j + \sum_{s=i'}^{n} p_{\pi'(s)} - \text{M}(\xi_{j})\\
    &\geq \bar{r}_{j^{*}} + \sum_{s=i^{*}}^{n} p_{\pi(s)} - \text{M}(\xi_{j^{*}})\\
    &= R(\pi, \xi_{j^{*}}) = R(\pi).
  \end{align*}
  The third inequality holds since permutation~$\pi$ is sorted in non-decreasing order of $\bar{r}_j - \text{M}(\xi_j)$ and $j$ is chosen as the smallest index in~$\pi'$. Since $j^{*}$ is a critical-regret job for~$\pi$ under scenario~$\xi_{j^{*}}$, the second equality is implied.\qed
\end{pf}

Again, the results of the previous theorem allow us to conclude polynomial time solvability of Robust Regret $1|r_j|C_{\max}$ with uncertainty set~$U_2$.  As in the previous section, it can be readily  verified that the result of Theorem~\ref{alg_regret} remains valid if we replace the uncertainty set~$U_2$ by~$U_1$.


In order to compute all values $M(\xi_j)$ we can first sort all the jobs in order of their lowest possible release dates in time $O(n\log n)$ and then update the sorted order $n$~times by raising the release date of a single job (and then decreasing it afterwards).  Each update of the order needs $O(\log n)$ time if we use binary search. Hence, all the optimum sequences~$M(\xi_j)$ can be computed in time~$O(n\log n)$. However, we still need to evaluate the objective function value for each of these solutions, which takes $O(n)$~time per solution, resulting in a total time of~$O(n^2)$. 

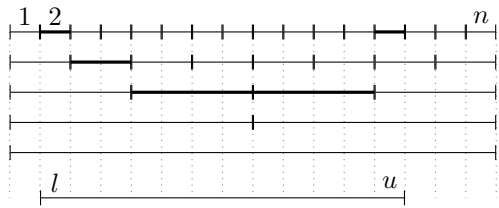
\begin{figure}
        \centering
        \begin{tikzpicture}[scale=0.8]
          \foreach \x in {0,...,16} {
            \draw[dotted,gray] (1+\x*0.5,0) -- (1+\x*0.5,-2.75);
          }
          \node at (1.25,0.25) {$1$};
          \node at (1.75,0.25) {$2$};
          \node at (8.75,0.25) {$n$};

          \foreach \x in {0,...,15} {
            \draw [|-|] (1+\x*0.5,0) --(1.5+\x*0.5,0);
          }

          \foreach \x in {0,...,7} {
            \draw [|-|] (1+\x,-0.5) --(2+\x,-0.5);
          }

          \foreach \x in {0,...,3} {
            \draw [|-|] (1+\x*2,-1) --(3+\x*2,-1);
          }

          \foreach \x in {0,...,1} {
            \draw [|-|] (1+\x*4,-1.5) --(5+\x*4,-1.5);
          }
          \draw [|-|] (1,-2) --(9,-2);

          \draw [|-|] (1.5,-2.75) --(7.5,-2.75);
          \node at (1.75, -2.5) {$l$};
          \node at (7.25, -2.5) {$u$};
          \draw[very thick] (1.5,0)--(2,0);
          \draw[very thick] (2,-0.5)--(3,-0.5);
          \draw[very thick] (3,-1)--(7,-1);
          \draw[very thick] (7,0)--(7.5,0);
        \end{tikzpicture}
        \caption{Illustration of the two-phase approach to compute the values~$m_j$.}
        \label{fig:construction}
      \end{figure}

We can, in fact, improve this running time to $O(n \log n)$ by describing a more efficient way of computing the values $M(\xi_j)$ for all jobs $j \in \mathcal{J}$, which is done in the lemma below. This enables us to solve the robust regret version in asymptotically the same time as the non-robust case.
\begin{lemma}
	The values $M(\xi_j)$ for all $j\in \mathcal{J}$ can be computed in $O(n\log n)$ time.
\end{lemma}
\begin{pf}
	We denote the scenario in which all release dates take their smallest value by \uxi, that is, $r_j^{\uxi} = \underline{r}_j$ for all $j\in\mathcal{J}$.
	First, we sort the jobs in non-decreasing order of these release dates in $O(n\log n)$ time and assume that $\underline{r}_1\leq\ldots\leq\underline{r}_n$ for the rest of the proof.
	In this situation we know that the optimal schedule for \uxi{} is given by the ERD rule and is the identity permutation $\pi$.
	This lets us compute information about the optimal schedule.
	We obtain, for all $j\in\mathcal{J}$, the completion times $C_j(\pi,\uxi) = \max\{C_{j-1}(\pi,\uxi),\underline{r}_j\}+p_j$ (where $C_0(\pi,\uxi) :=0$) and the differences $\varepsilon_j:= C_j(\pi,\uxi) - (\underline{r}_j + p_j)$ between these and the earliest possible completion times.
	Additionally we compute the idle time $a_j = \max\{\underline{r}_j-C_{j-1}(\pi,\uxi) ,0\}$ directly before job $j$ and the total idle time $f_j = \sum_{s>j}a_s = f_{j+1} + a_{j+1}$ after job $j$, for all $j\in\mathcal{J}$.
	All of this requires linear time.
	
	For $j\in \mathcal{J}$, let $u_j:=\max\{ k : \underline{r}_k \leq \overline{r}_j \}$, which can be computed in $O(\log n)$ time for any fixed $j$ by a binary search.
	Let $\pi_j$ be the permutation where we move the job $j$ to the position after $u_j$, that is, $\pi_j(i) = i$ for $i\in \{1,\ldots,j-1, u_j+1,\ldots,n\}$, $\pi_j(i) = i+1$ for $i\in \{j,\ldots,u_j-1\}$, and $\pi(u_j) = j$.
	This permutation is optimal for $\xi_j$ by the ERD rule.
	One final parameter we require is $m_j := \min \{p_j,\varepsilon_{j+1},\ldots,\varepsilon_{u_j}\}$.
	We show how to compute these values efficiently at the end of the proof. 
	
	To determine $M(\xi_j)$ we note that $C_s(\pi_j,\xi_j) = C_s(\pi,\uxi)$ for $s=1,\ldots, j-1$.
	The removal of job $j$ at the $j$th position means that we can potentially schedule the jobs $s$ in the interval $j+1,\ldots,u_j$ sooner.
	Job $j+1$ can be scheduled earlier by $\min\{p_j, \varepsilon_{j+1}\}$ units of time and subsequent jobs $s$ are limited by $\varepsilon_s$ and the decrease of the previous job's completion time.
	Inductively this gives us that $C_{s}(\pi_j,\xi_j) = C_{s}(\pi, \uxi) - \min\{p_j,\varepsilon_{j+1},\ldots,\varepsilon_s\}$.
	Consequently, the completion time of job $j$ under $\xi_j$ is given by
	\begin{align*}
	C_j(\pi_j,\xi_j) 
	&= \max \{ C_{u_j}(\pi_j,\xi_j), \overline{r}_j \} + p_j
	= \max \{ C_{u_j}(\pi,\uxi) - m_j, \overline{r}_j \} + p_j \\
	&= C_{u_j}(\pi,\uxi) + p_j - m_j + \max\{ 0, \overline{r}_j - C_{u_j}(\pi,\uxi) + m_j\} 
	=: C_{u_j}(\pi,\uxi) + h_j.
	\end{align*}
	The remaining jobs $u_j +1,\ldots,n$ are scheduled in the same order as before, with the sole difference that the job $u_j$ is now delayed by $h_j\geq 0$.
	We note that the delay $z_s := C_s(\pi_j, \xi_j) - C_s(\pi,\uxi)$ of job $s \in\{ u_j+1,\ldots,n\}$ is then equal to the delay of the job preceding it, decreased by the potential idle time $a_s$ before processing job~$s$.
	Therefore $z_{u_j+1} = \max \{h_j-a_{u_j+1},0\}$ and $z_s = \max \{z_{s-1} - a_s, 0\}$ for $s>u_j+1$.
	The last expression can be rewritten to $z_s = \max \{h_j -\sum_{k=u_j+1}^{s}a_j,0\}$ resulting in $M(\xi_j) = M(\uxi) + \max\{h_j - f_{u_j}, 0\}$.
	Using this formula and our previously computed parameters we can determine each $M(\xi_j)$ in constant time.
	
	So we just need to show how to compute the values $m_j$ in $O(n\log n)$ time.
	We describe a pre-processing that requires $O(n)$ time and then give an algorithm to determine $\min \{\varepsilon_s : l\leq s \leq u \}$ in $O(\log n)$ time.
	This lets us compute $m_j$ as the minimum of $p_j$ and $\min\{\varepsilon_s : j+1 \leq s \leq u_j \}$. 
	
	In the preprocessing step we compute the minima in intervals of length $2, 4, \ldots, 2^q$ where $2^{q-1} < n \leq 2^q$.
	More precisely, in step $k$ we compute the minima in the intervals from $j+1$ to $j+2^k$ where $j$ satisfies that $2^k$ divides $j$ ($2^k\mathrel{|}j$) and $j+2^k\leq u$.  The intervals are illustrated in Figure~\ref{fig:construction}.
	These computations require constant time as  it is just the minimum of the two minima of the intervals $j+1,\ldots,j+2^{k-1}$ and $j+2^{k-1} + 1,\ldots,j+2^k$ computed in the previous step.
	As there are $\lfloor \tfrac{n}{2^k}\rfloor$ intervals of length $2^k$ that we need to compute, the total running time is $O(n)$.
	
	The algorithm for computing $\min\{\varepsilon_s : l \leq s \leq u \}$ now uses the precomputed values by, essentially, taking the minimum over largest interval available starting at the lower bound that does not exceed $u$ and then updating the lower bound to the end of that interval (plus one).
	It then returns the minimum of all values regarded, which is the result.
	
	To obtain the running time, we run the algorithm in two phases.
	The first is an increasing phase, in which the lengths of the intervals regarded are non-decreasing.
	In it, we start at $j = l-1$ by regarding intervals from $j+1$ to $j+2^k$ for $k=0$.
	As long as $j+2^{k+1}\leq u$ and the value for this longer interval is precomputed, so if $2^{k+1}\mathrel{|} j$, we increment $k$.
	If one of these conditions fails, we store the minimum on the interval $j+1,\ldots,j+2^k$ and increase $j$ by $2^k$ and continue.
	Note that $2^k \mathrel{|} j$ implies that $2^k$ also divides $j+2^k$, so the interval of length $2^k$ is also precomputed at the new $j$.
	This phase ends when $j + 2^k > u$, at which point we require smaller intervals to cover the remainder.
	
	This phase needs $O(\log n)$ time as $k$ is bounded from above by $q\leq \log n$ and there can be at most that many increment steps.
	An update step caused by $j+2^{k+1} > u$ results in a termination of the phase, whereas one caused by $j$ not being divisible by $2^{k+1}$ results in the new index $j+2^{k}$ being a multiple of $2^{k+1}$.
	This means that such an update is followed by an increment step or a violation of $j+2^{k+1} \leq u$ in the next step, the latter of which leads to a termination once more.
	
	In the second, decreasing, phase we decrement $k$ if $j+2^k > u$ and otherwise perform an update by again storing the minimum from $j+1$ to $j+2^k$ and increasing $j$ by $2^k$.
	Here every update reduces the length of the remaining interval by at least half, and there are are most $O(\log n)$ updates and decrements.
	Hence, this procedure computes the minimum in logarithmic time.\qed
      \end{pf}

\begin{corollary}
	The problem Robust Regret $1|r_j|C_{\max}$ with uncertainty set~$U_1$ or~$U_2$ can be solved optimally in time~$O(n\log n)$.
\end{corollary}



\end{document}